\documentclass{article}

\begin{document}

\title{Chicken or egg? A hierarchy of homotopy algebras}         
\author{F\"{u}sun Akman}        
\date{\today}          
\maketitle

\begin{abstract} 
We start by clarifying and extending the multibraces notation, which economically describes substitutions of multilinear maps and tensor products of vectors. We give definitions and examples of weak homotopy algebras, homotopy Gerstenhaber and Gerstenhaber bracket algebras, and homotopy Batalin-Vilkovisky algebras. We show that a homotopy algebra structure on a vector space can be lifted to its Hochschild complex, and also suggest an induction method to generate some of the explicit (weakly) homotopy Gerstenhaber algebra maps on a topological vertex operator algebra (TVOA), their existence having been indicated by Kimura, Voronov, and Zuckerman in 1996 (later amended by Voronov). The contention that this is the fundamental structure on a TVOA is substantiated by providing an annotated dictionary of weakly homotopy BV algebra maps and identities found by Lian and Zuckerman in 1993.
\end{abstract}

\tableofcontents

\section{Introduction}       

Following the work of Kimura, Voronov, and Zuckerman \cite{KVZ} in 1996, we were able to reduce the defining identities of a certain homotopy Gerstenhaber algebra structure on a topological vertex operator algebra (TVOA) to a statement about the composition of mega operators. Later, Voronov \cite{Vor} made some corrections to \cite{KVZ}, and as a result, our construction was only partially applicable (although the same mega identity continues to fit perfectly the structure on the Hochschild complex of an associative algebra).  

We will continue to regard the construction in \cite{KVZ} and \cite{Vor} as a weakly homotopy Gerstenhaber algebra. Unless specifically designated ``strongly'' homotopy algebras (an operadic, or according to M.~Markl, ``God-given'', property), all homotopy algebras mentioned in the article are ``weakly'' homotopy algebras. That is, they have the algebraic property that some of the constituent maps descend to a classical structure in the cohomology. Moreover, ``Gerstenhaber'' and ``Batalin-Vilkovisky'' will be shortened to G and BV respectively, where appropriate.

We are now in a position to provide a characterization for the intriguing homotopy Batalin-Vilkovisky algebra structure on a TVOA, discovered several years earlier by Lian and Zuckerman \cite{LZ} and published in 1993. True to the chronological order, it has been conjectured that the anti-symmetrization of the homotopy BV maps of Lian and Zuckerman will lay some fresh homotopy Gerstenhaber maps \cite{KVZ}. We refer the reader to Huang and Zhao's article \cite{HZ} on a new operadic and geometric formulation of a TVOA and the proof of the conjecture by Lian and Zuckerman and Kimura et al.~stating that a TVOA gives rise to a homotopy Gerstenhaber algebra. It turns out that the egg is the parent of the chicken after all; we can identify some components of the lower identities in \cite{LZ} as bona fide homotopy G-operators in the sense of \cite{KVZ} and \cite{Vor}. In fact, every explicit map introduced by Lian and Zuckerman as part of the homotopy BV structure is simply the bracket of the anti-ghost operator $b_0$ with the homotopy G-maps, and every explicit identity is either homotopy Gerstenhaber, or can be obtained from one by bracketing with $b_0$. Then why not define a homotopy BV algebra to be in general an algebra obtained by bracketing homotopy G-maps with a suitable odd operator? (See Tamarkin and Tsygan's definitive work on $BV_{\infty}$ algebras \cite{TT}.) The correspondence is visible in \cite{LZ} once each identity is deciphered as ``sums of compositions of maps equals zero''. Lian and Zuckerman's work also gives us an idea as to how to construct some of the homotopy G-algebra maps on a TVOA explicitly by induction. After establishing a reasonable definition of homotopy BV algebras, we show that any homotopy algebra structure on a space $V$ can be lifted to its Hochschild complex, $C^{\bullet}(V,V)=\mbox{\it \mbox{\it Hom}}(TV,V)$ (Getzler has a similar construction for $A_{\infty}$ algebras in \cite{Get}).

The present paper originated from the author's confusion about the proliferation of G-objects, especially on a TVOA, and a desire to develop and clarify the multibraces language (\cite{Akm2}, \cite{Akm3}) to the point that it effortlessly describes and advances most homotopy algebra structures so far discovered. Pairs of braces $\{\; \}\{ \; \}$, as well as $\{\; \}\{ \; ,\; ,\dots ,\; \}$ were introduced by Gerstenhaber \cite{Ger}, and much later revived in more general form  by Kadeishvili and Getzler. Our multibraces are a generalization of this idea, allowing one to substitute any number of maps into each other. Moreover, by defining the symbol $\{ a_1,\dots, a_n\}$ to be the tensor product of the vectors $a_1,\dots,a_n$ (as it should be), we extend the idea of substitutions of maps from a tensor product $TV$ to the vector space $V$, and of vectors into these maps, to substitutions and tensor products in $\mbox{\it Hom}(TV,TV)$. Whenever multilinear maps are defined on $\mbox{\it Hom}(TV,V)$ (or $\mbox{\it Hom}(TV,TV)$) itself, we use the ``second-level braces'' \cite{Akm2} with primes. 

Throughout the paper we point to the dichotomy between mega maps that are called ``(weakly) homotopy'' because the structure descends to the classical case in the cohomology, and mega maps that should be called ``homotopy'' because they satisfy identities similar to the classical case (and descend to it in the cohomology as well). 

\section{Substitution operators}

\subsection{Multibraces}

\subsubsection{Substitutions of multilinear maps}    

We will not make any references to operads in our description; the exposition by Jones in \cite{Jon} is recommended for a very refined treatment of substitution operators, and much more, in the context of operads. All our algebraic structures will be defined on a real or complex vector space $V$ with one or more gradings to be discussed later. Let us agree to denote specific $n$-linear maps from $V\times\cdots\times V$ into $V$ by $(n)$ and in particular vectors in $V$ by $(0)$. The symbol
\[ (\, (n_1),(n_2),\dots,(n_k)\, )\]
denotes the substitution operator that places specific $n_1$-linear, $n_2$-linear, ..., $n_k$-linear maps (possibly just vectors) into a $k$-linear map ($k\geq 1$) in this order. The resulting composite map is then $\sum_{t=1}^kn_t$-linear.

It is also desirable to leave certain arguments of the outer, $k$-linear map empty while substituting specific maps and vectors into others. If $r$ consecutive arguments are thus selected, we denote this by the unadorned number $r$ instead of $(r)$. The most general substitution operator creates a map whose arguments can be counted by adding up all the integers that appear inside its symbol, with or without parentheses.

{\bf Example 1.} The symbol $(\, (2),(1),4,(3),(0)\, )$ denotes the substitution of a bilinear map and a linear map into the first two arguments of an 8-linear map, followed by four blank spaces, a trilinear map, and a specific vector. The final map has room for $2+1+4+3+0=10$ vectors.

The composition of two substitution operators has a simple rule. First consider two operators $S_1=(\, (n_1),\dots,(n_k)\, )$ and $S_2=(\, (m_1),\dots,(m_l)\, )$, where all entries are within parentheses, and $n_1+\cdots +n_k=l$. Their composition is given by  
\[ S_1\circ S_2=\left( \, \left( \sum_{i=1}^{n_1}m_i\right),\left( \sum_{i=1}^{n_2}m_{n_1+i}\right),\dots, \left( \sum_{i=1}^{n_k}m_{n_1+\cdots +n_{k-1}+i}\right) \, \right) ,\]
that is, we add up the first $n_1$ entries on the right, then the next $n_2$, and so on. When $n_i=0$ we have an empty sum and the resulting entry is $(0)$. 

{\bf Example 2.} $(\, (2),(1),(3),(0)\, ) \,\circ\, (\, (1), (0), (3), (0), (2), (4)\, ) = (\, (1),(3),(6),(0)\, )$.

In case the operator on the left involves positive integers $r$ without parentheses, these numbers are to be copied as is, since they correspond to empty arguments. 

{\bf Example 3.} $(\, (2),(1),3\, ) \,\circ\, (\, (1), (5), (4)\, ) = (\, (6),(4),3\, )$.

Similar numbers $r\geq 1$ on the right should be converted to $r$ 1's separated by commas, and then the usual rules are applied, as if all entries are surrounded by parentheses.

{\bf Example 4.} $(\, (1), (3)\, ) \,\circ\, (\, (2),(4),2\, )=(\, (1), (3)\, ) \,\circ\, (\, (2),(4),1,1\, ) =(\, (2),(6)\, )$.

\subsubsection{Bigrading in the Hochschild space}

Let 
\[ V=\oplus_{j\in\mbox{\bf Z}}V^j\]
be a {\bf Z}-graded vector space over the complex field. We will call this grading ``super'', and assume that any multilinear map under consideration is a finite sum of homogeneous maps under the super grading throughout the paper. The notation for super degree will be 
\[ \mbox{$|a|=j$ if $a\in V^j$ and $|m|=j$ if $|m(a_1,\dots,a_n)|=|a_1|+\cdots +|a_n|+j$} \]
for all homogeneous $a_i\in V$. The expressions ``odd operator'' or ``even operator'' will refer to the parity of a homogeneous map under super grading. Another simple notion of degree, namely the number of arguments of a multilinear map minus one, will be given by
\[ \mbox{$d(m)=n-1$ if $m:V^{\otimes n}\rightarrow V$}.\]
Note that vectors, as 0-linear maps, have $d$-degree negative one.

An $n$-linear map $m$ will be called {\it symmetric} if 
\[ m(a_1,\dots, a_i,a_{i+1},\dots,a_n)=m(a_1,\dots,a_{i+1}, a_i,\dots,a_n),\]
{\it anti-symmetric} if
\[ m(a_1,\dots, a_i,a_{i+1},\dots,a_n)=-m(a_1,\dots,a_{i+1}, a_i,\dots,a_n),\]
{\it super symmetric} if
\[ m(a_1,\dots, a_i,a_{i+1},\dots,a_n)=(-1)^{|a_i|\, |a_{i+1}|}m(a_1,\dots,a_{i+1}, a_i,\dots,a_n),\]
and {\it super anti-symmetric} if
\[ m(a_1,\dots, a_i,a_{i+1},\dots,a_n)=-(-1)^{|a_i|\, |a_{i+1}|}m(a_1,\dots,a_{i+1}, a_i,\dots,a_n)\]
for $1\leq i\leq n-1$. If there is more than one grading on the vector space, say $|\; |_1$ and $|\; |_2$, then bigraded symmetry and bigraded anti-symmetry are defined by replacing $|a_i|\, |a_{i+1}|$ by $|a_i|_1\, |a_{i+1}|_1+|a_i|_2\, |a_{i+1}|_2$. 

In case of a generic vector space $V$, we call the space
\[ C^{\bullet}(V,V)=\mbox{\it Hom}_{\mbox{\bf C}}\left( \bigoplus_{n=0}^{\infty}  V^{\otimes n};V\right)  \]
of multilinear maps on $V$ the {\it Hochschild space}. If there is additional structure and/or a distinguished differential on $C^{\bullet}(V,V)$, then the names {\it Hochschild algebra} and {\it Hochschild complex} will be used respectively. The bigrading on $C^{\bullet}(V,V)$ is given by the super degree and the $d$-degree defined above. 

The space $C^{\bullet}(V,V)$ consists of maps $m$ which can be decomposed into infinitely many $d$-homogeneous parts $m_{(n)}$ by restriction to $T^nV$. Later, we will consider the {\it partitioned Hochschild space} with maps defined on the {\it partitioned tensor space} $T_{par}V$, where, for example, we will distinguish between basic vectors $a\otimes b\otimes c$, $(a)\otimes(b\otimes c)$, ... because of the way the tensor factors are grouped. We will denote the corresponding direct summands of $T_{par}V$ by $T^{(3)}V$, $T^{(1|2)}V$, and so on. Then maps that reside inside the partitioned Hochschild space will be formal sums of $m_{\pi}$ where $\pi$ is a partition of some positive integer.

\subsubsection{Substitutions with braces}

The expression \cite{Akm2}
\[ \{ x\} \{ y_1,\dots,y_k\} \cdots \{ z_1,\dots z_l\}\{ a_1,\dots,a_n\} \]
where the $a_i$'s are vectors and the remaining letters denote specific multilinear maps, possibly vectors, means that every symbol except $x$ is to be substituted into one on the left, in every possible way, but two vectors as well as two symbols within the same pair of parentheses may not to be combined in this manner. The resulting substitutions may not change the order within any one pair of parentheses, either. The signed sum of all such expressions are then formed, taking into consideration the bidegrees of the symbols that have been interchanged with respect to the original above. Note that although we require there be enough arguments at every level to accommodate all symbols to the right, we do {\it not} require an exact match. In other words, there may be some free arguments (to the left of the string of $a_i$'s) that result in our sliding the multilinear maps within the boundaries of the pecking order. In this sense, braces denote sums of the inflexible substitutions described in the previous section, although the leftmost map $x$ was only implicitly represented, and more than two nested maps were not considered.

Also note that the $a_i$'s are not an essential part of braces, just as the composition $f\circ g$ of ordinary functions is meaningful on its own without an argument. If included, they are frequently used to display the number (and grouping) of arguments of the composite multilinear function to their left. Let us establish here the convention of reserving the symbols $a_i$ (or $a,b,c,\dots$) for this exclusive purpose.

We can now define the plus/minus sign picked up in the exchange of maps (possibly vectors) $x$ and $y$ to be
\[ (-1)^{|x|\,|y|+d(x)d(y)}.\]
We also define the {\it Gerstenhaber bracket} of two maps $x$ and $y$ by
\[ [x,y]=\{ x\}\{ y\} -(-1)^{|x|\,|y|+d(x)d(y)}\{ y\}\{ x\} ,\]
as Gerstenhaber did in \cite{Ger}.

{\bf Example 1.} Given a bilinear map $x$ and linear and four-linear maps $y$, $z$ respectively, the substitution $\{ x\}\{ y,z\}$ is of type $(\, (1),(4)\, )$. On the other hand, $\{ x\}\{ z\}$ gives us two choices for the placement of $z$, thus we obtain substitutions of types $(\, (4),1\, )$ and $(\, 1, (4)\, )$. Written with arguments, the substitutions mentioned above are given by
\begin{eqnarray} && \{ x\}\{ y ,z\}\{ a_1,a_2,a_3,a_4,a_5\}  \nonumber \\ &=&(-1)^{|a_1|\, |z|+d(a_1)d(z)} x(\, y (a_1),z(a_2,a_3,a_4,a_5)\, ) \nonumber \\
&=& -(-1)^{|a_1|\, |z|}x(\, y (a_1),z(a_2,a_3,a_4,a_5)\, ) \nonumber \end{eqnarray}
and
\begin{eqnarray} && \{ x\}\{ z\}\{  a_1,a_2,a_3,a_4,a_5\} \nonumber \\ &=& x(\, z(a_1,a_2,a_3,a_4),a_5\, )\nonumber \\ & &-(-1)^{|a_1|\, |z|} x(\, a_1,z(a_2,a_3,a_4,a_5)\, )\nonumber \end{eqnarray}
respectively.

{\bf Example 2.} $\{ x\}\{ y\}\{ z\}$ means that $y$ is substituted into $x$ in all possible ways and $z$ may either go into $y$ or straight into $x$. That is, we have
\[ \{ x\}\{ y\}\{ z\} =\{ x\}\{ \{ y\}\{ z\}\} +\{ x\}\{ y,z\} \pm \{ x\}\{ z,y\} .\]

{\bf Example 3.} Let $m$ be an even bilinear map. The identity $\{ m\}\{ m\} =0$ is nothing but the associativity of $m$:
\[ \{ m\}\{ m\}\{ a_1,a_2,a_3\} =m(\, m(a_1,a_2),a_3\, )-m(\, a_1,m(a_2,a_3)\, )=0.\]
When the name of the map is suppressed, we obtain the familiar statement $(a_1a_2)a_3=a_1(a_2a_3)$.

{\bf Example 4.} By definition, $\{ m\}\{ a_1\}\{ a_2\}\{ a_3\}$ denotes the anti-symmetrization of the expression $\{ m\}\{ a_1,a_2,a_3\}$ for a trilinear map $m$:
\[ \{ m\}\{ a_1\}\{ a_2\}\{ a_3\} =\sum_{\sigma\in S_3}\mbox{sgn}(\sigma) m(a_{\sigma(1)},a_{\sigma(2)},a_{\sigma(3)}),\]
where $S_3$ denotes the symmetric group on three letters. If there is a super degree, a factor of $(-1)^{|a_i|\, |a_j|}$ is to be inserted for every crossing of symbols $a_i$, $a_j$. 

We observe first that both the super degree and the $d$-degree are preserved by braces. Second, the braces $\{ \, \}\{ \, \}$ form a right pre-Lie product on the Hochschild space. A bilinear product $\star$ is called {\it right pre-Lie} if
\[ (a\star b)\star c -a\star (b\star c)=(-1)^{|b|\, |c|}((a\star c)\star b-a\star (c\star b))\]
holds for all homogeneous $a,b,c$. Gerstenhaber proved this result in detail in \cite{Ger}, but see \cite{Akm2} for a very brief proof based on the new braces notation. As a result, the Gerstenhaber bracket makes the Hochschild complex into a bigraded Lie algebra.

\subsubsection{Isomorphism between Hochschild space and coderivations on the tensor coalgebra}\label{isom}

The use of multibraces, being very visual, simplifies the statements and proofs of many results. For future reference, and to set the notation, let us look into the isomorphism in the section heading. 

Obviously, $\{ a_1,\dots,a_n\}$ means the tensor product $a_1\otimes\cdots a_n\in V^{\otimes n}\subset TV$. This expression, so far, lives only inside a multilinear map. Also expressions of type
\[ \{ m_{(2)}\} \{ a,b,c\} =\{  m_{(2)}(a,b),c\} \pm \{ a,  m_{(2)}(b,c)\} \in V^{\otimes 2}\subset TV, \]
where the maps are ``smaller'' than the arguments, are undefined in our current setting on their own, although this latter would make perfect sense inside some $\{ m_{(3)}\}$. If elements of $\mbox{\it Hom}(TV,TV)$ are allowed, though, we will have a lot of use for braces. Their use is completely symmetric between multilinear maps and tensor products of vectors. Note that terms such as
\[ \{ m_{(2)}(a,-),b,c\}\]
are not permitted in these constructions; that is, vectors are required to fill consecutively all possible slots inside those multilinear maps to their left and in another pair of braces.

The {\it second-level braces}, namely those denoting tensor products and substitutions in the tensor algebra $T(TV)$, were also introduced in \cite{Akm2}. We distinguish this level by putting primes on the braces and the tensor product symbols. (In fact, we will use primed objects whenever multilinear maps are defined on the Hochschild complex itself.) Since we want to study coderivations and the comultiplication, it makes sense to consider maps from $TV=T^1(TV)$ into $T(TV)$. Let us now denote the tensor product $a_1\otimes\cdots\otimes a_n$ inside $T^1(TV)$ by
\[ \{ \, \{ a_1,\dots,a_n\} \, \}' .\]
For example, the basis element
\[ (a_1\otimes \cdots \otimes a_n)\otimes'(b_1 \otimes \cdots \otimes b_k)          \]
of $T^2(TV)$ will be written as
\[ \{ \, \{ a_1,\dots,a_n\} , \{ b_1,\dots,b_k\}\,\}' .    \]
Similarly, an old multilinear map $m$ will be denoted by $\{ \, \{ m \} \, \}'$, and substitutions between symbols of the first level will still take place even if they are in different pairs of primed braces: we have
\[ \{ \, \{ m \} \, \}' \{ \, \{ a_1,\dots,a_n\} \, \}' = \{ \{ m\}\{ a_1,\dots,a_n\} \, \}'\]
by definition. One such map we will use is the identity map
\[ i=i_{(0)}+i_{(1)}+i_{(2)}+\cdots +i_{(n)} +\cdots :TV\rightarrow TV \]
where $i_{(n)}( a_1\otimes \cdots \otimes a_n)=a_1\otimes \cdots \otimes a_n$. ``Primed'' maps will be given capital letter names to distinguish them from the first-level maps. 

The standard associative bilinear product $M$ on $TV$,
\[ M: T^1(TV)\otimes T^1(TV)\rightarrow T^1(TV),\]
 is defined via $i$:
\begin{eqnarray} && \{ M\}' \{ \, \{ a_1,\dots,a_n\} , \{ b_1,\dots,b_k\}\,\}' \nonumber \\
&=& \{\,\{   i\}\,  \}' \{ \, \{ a_1,\dots,a_n\} , \{ b_1,\dots,b_k\}\,\}' \nonumber \\
 &=& \{ \, \{ i_{(n+k)}\}\, \}' \{ \, \{ a_1,\dots,a_n\} , \{ b_1,\dots,b_k\}\,\}' \nonumber \\
&=& \{ \, \{ a_1,\dots,a_n ,  b_1,\dots,b_k\}\,\}' .\nonumber \end{eqnarray}
Meanwhile, the coassociative coproduct $\Delta$ on $TV$,
\[ \Delta : T^1(TV)\rightarrow T^1(TV)\otimes T^1(TV),\]
has a similar definition:
\begin{eqnarray} && \{ \Delta\}' \{ \, \{ a_1,\dots,a_n\} \,\}' \nonumber \\
&=& \{ \, \{ i,  i\}\,  \}' \{ \, \{ a_1,\dots,a_n\} \,\}' \nonumber \\
 &=&  \sum_{u+v=n} \{ \, \{ i_{(u)},i_{(v)}   \}\, \}' \{ \, \{ a_1,\dots,a_n\} \,\}' \nonumber \\
&=&  \sum_{u=0}^{n} \{ \, \{ a_1,\dots,a_u\} , \{ a_{u+1},\dots,a_n\}\,\}' .\nonumber \end{eqnarray}
Then in accordance with the usual rules, we define  {\it derivations}  and  {\it coderivations}  of the bialgebra $TV$ as 
\[ \mbox{\it Der}(TV)=\{ X\in \mbox{\it Hom}(T^1(TV),T^1(TV)): [X,M]'=0\}\]
and
\[ \mbox{\it Coder}(TV)=\{ X\in \mbox{\it Hom}(T^1(TV),T^1(TV)): [X,\Delta ]'=0\} .\]

It is now easy to see why $C^{\bullet}(V,V)$ and $\mbox{\it Coder}(TV)$ are naturally isomorphic. Given 
\[ m=m_{(1)}+m_{(2)}+\cdots m_{(n)}+\cdots \in \mbox{\it Hom}(TV,V),\]
define $\delta (m)\in \mbox{\it Coder}(TV)$ by
\begin{eqnarray} && \{ \delta (m) \}' \{ \, \{ a_1,\dots,a_n\} \,\}' \nonumber \\
 &=& \{ \, \{ m \} \, \}' \{ \, \{ a_1,\dots,a_n\} \, \}' \nonumber \\
&=& \sum_{k=1}^n \{ \, \{ m_{(k)} \} \, \}' \{ \, \{ a_1,\dots,a_n\} \, \}' \nonumber \\
&=& \sum_{k=1}^n \sum_{j=1}^{n-k+1}\pm \{ \, \{ a_1,\dots, m_{(k)}(a_j,\dots, a_{j+k-1}),\dots, a_n \} \, \}' .\nonumber \end{eqnarray}
Why is $\delta (m)$ a coderivation? We have
\begin{eqnarray} &&   \{ \Delta \}' \{ \delta (m) \}'\{ \, \{ a_1,\dots,a_n\} \, \}'\nonumber \\
&=& \{ \, \{ i,  i\}\,  \}' \{ \, \{ m\} \, \}'  \{ \, \{ a_1,\dots,a_n\} \,\}' .\nonumber \end{eqnarray}
and 
\begin{eqnarray} &&  \{ \delta (m) \}' \{ \Delta \}' \{ \, \{ a_1,\dots,a_n\} \, \}'\nonumber \\
&=& \{ \, \{ m\} \, \}' \{ \, \{ i,  i\}\,  \}' \{ \, \{ a_1,\dots,a_n\} \,\}'\nonumber \end{eqnarray}
In the first case, the first $n$ summands of the mega operator $m$ are applied in all possible ways to $\{ a_1,\dots,a_n\}$, then the resulting product is split into two in all possible ways as tensor products in $T^2(TV)$. In the second case, first $\{ a_1,\dots,a_n\}$ is split into two in all possible ways as tensor products in $T^2(TV)$, then the same summands of $m$ are applied in all possible ways to the split product. The end result is the same!

Going backwards from $\mbox{\it Coder}(TV)$, we take any coderivation and compose it with the projection $TV\rightarrow V$, which gives a map in $\mbox{\it Hom}(TV,V)$. Note that coderivations on $TV$ are uniquely determined by this projection, just as derivations would be uniquely determined by the restriction to $V$.

Moreover, the reason why Gerstenhaber's bracket on $\mbox{\it Hom}(TV,V)$ is the same as the bracket of coderivations as linear maps is the following: both are written as $\{ x\}\{ y\} \pm \{ y\}\{ x\}$, and are subject to the multibraces rules. 

\subsection{Partitioned multibraces}

\subsubsection{Substitutions of partitioned multilinear maps} \label{parmaps}   

We are back to the round-bracket substitutions, where every single map and vector is placed just so. We will now replace the notation $(n)$ denoting the degree of linearity of a multilinear map by a partition $\pi =(j_1|j_2|\cdots |j_s)$ ($s\geq 1$), where the ``parts'', or ``slots'', $j_1$, $j_2$, ..., $j_s$ add up to $n$. This change affects only the composition properties of the function. The simplest case to consider is $k$ {\it partitioned maps} substituted into an ordinary $k$-linear map. We write, as before, $(\pi_1,\pi_2,\dots, \pi_k)$, and interpret the result as the partitioned map whose type is obtained by ``merging'' all partitions. In other words, we erase all inner parentheses, but leave the bars intact, and replace commas by plus signs.

{\bf Example 1.} The result of the substitution $(\, (1|2),(0|2|5),(9)\, )$ is a $(\, 1|2|2|14 \, )$-linear map.

If partitioned maps are substituted into a partitioned map instead of a plain one, the procedure is the same. Bars at both levels are left in place and addition is performed only at the commas.

{\bf Example 2.} The substitution $(\, (1|2),(0|2|5)|(9), (3|0|4), (1|1)| (3|2|2)\, )$ produces a $(\, 1|2|2|5|12|0|5|1|3|2|2\, )$-linear map.

What happens if there are ``free'' arguments of the outer map at which there is no substitution? What is the type of the new map? Generalizing from the plain-map case, we conclude that these entries need to be treated as if they are enclosed by parentheses.

{\bf Example 3.} The substitution $(\, (1|2),3,(5|7)|6\, )$ gives us a $(\, 1|10|7|6\, )$-linear map.

So far, the reason for finer partitioning among the arguments of a multilinear map is not apparent. When barred braces are used for substitution, though, the types of multilinear functions that are produced play a prominent role in writing homotopy algebra identities.

Partitions with 0's will be used to write the symmetries of a map. For example, a multilinear map of type $(1|0)$ is in fact a linear map that only acts on the first slot and ignores the slot immediately to the right.   

\subsubsection{Substitutions with partitioned braces} \label{parbraces}

The algebra of partitions and partitioned braces first came up in the author's work describing a master identity for homotopy Gerstenhaber algebras \cite{Akm3} (later shown to be inadequate except for lower identities as a result of Voronov's amendment \cite{Vor} to \cite{KVZ}). They were the algebraic equivalent of the construction of boundary (differential) in the pictures of Fox-Neuwirth cells in \cite{KVZ} and \cite{Vor}. However, the master identity holds perfectly for the Hochschild complex of an associative algebra, as we shall see.

How do we compose a partitioned map (on the left) with another, or several others (on the right)? As long as we do not consider the composite map's evaluation patterns, the answer is similar to the one we would give for plain map compositions... place all maps on the right, necessarily within the same pair of braces, into the one on the left at the same time, in all possible ways. Be sure to preserve the order. Then add up the different configurations.

In order to understand the mechanics of the substitution
\[ \{ m_{\pi}\} \{ m_{\pi_1},\dots,m_{\pi_l}\} \{ a_1^{(1)}\cdots a_{n_1}^{(1)} | \cdots |  a_1^{(k)}\cdots a_{n_k}^{(k)} \} ,\]
consider the simplest case of
\[ \{ m_{(n)}\} \{ m_{(i_1 | \cdots | i_k)} \} \{ a_1^{(1)}\cdots a_{n_1}^{(1)} | \cdots |  a_1^{(k)}\cdots a_{n_k}^{(k)} \} . \]
That is, only one partitioned map is to be substituted into $m_{(n)}$ (hence $l=1$), and the number of slots in the partition defining the substitution properties of this inner map matches that of the partitioning of the vectors, namely $k$. If each slot of $ m_{(i_1 | \cdots | i_k)}$ is at most equal to the number of vectors in the corresponding slot of the ordered vector set, we are in business (otherwise substitution is not possible). In the notation of the last Eq., we want $i_j\leq n_j$ for all $1\leq j\leq k$. 

{\bf Example 1.a.} Let us start working on an example that we will keep coming back to throughout this section. Consider the substitution
\[ \{ m_{(4) }\} \{ m_{(1|3)}\}\{ a_1,a_2,a_3 | b_1,b_2,b_3,b_4\} .\]
Both the middle map and the set of arguments sport two slots ($k=2$), and in each vector slot there are more vectors than $m_{(1|3)}$ can hold. 

In the general case, we want to split each slot of vectors into three pieces, some of which may be $\emptyset$, in all possible ways. The only condition is that the middle piece be as long as the number in the corresponding slot for $\{ m_{(i_1 | \cdots | i_k)} \}$. That is, we are getting ready to plug in the middle set in each slot into its designated place inside $\{ m_{(i_1 | \cdots | i_k)} \}$. We may denote this splitting process by
\begin{eqnarray} \{ a^{(1)}\} &=& \{ b^{(1)},c^{(1)},d^{(1)}\} \nonumber \\
&\vdots & \nonumber \\
\{ a^{(k)}\} &=& \{ b^{(k)},c^{(k)},d^{(k)}\} , \nonumber \end{eqnarray}
where letters with superscripts only denote blocks of consecutive vectors. In particular, the $c$'s have the exact number of vectors to fit in the corresponding slot of $ m_{(i_1 | \cdots | i_k)} $.

{\bf Example 1.b.} The first slot of vectors in Example 1.a can be split as
\[ \emptyset \;\;\; \{ a_1 \} \;\;\; \{ a_2,a_3\} \;\;\; \mbox{or}\;\;\;  \{ a_1\}\;\;\; \{ a_2\} \;\;\; \{ a_3\} \;\;\; \mbox{or}\;\;\;  \{ a_1,a_2\} \;\;\; \{ a_3\} \;\;\; \emptyset    .\]
Similarly, the second slot offers the possibilities 
\[ \emptyset \;\;\;\{ b_1,b_2,b_3\} \;\;\; \{ b_4\}\;\;\; \mbox{and} \;\;\; \{ b_1\}\;\;\; \{ b_2,b_3,b_4\}\;\;\; \emptyset .\]
Therefore, there are altogether six different ways (splittings) in which we may choose to apply $m_{(1|3)}$ to the partitioned vector set. 

The next step is, for each overall choice of splittings, to permute the ``left leftovers'' in all slots of vectors in such a way that inner orders are preserved. In the same manner we permute all ``right leftovers'' so that {\it their} inner orders are preserved, and apply $\{ m_{(i_1 | \cdots | i_k)} \}$ to the $c$'s in such a way that the left leftovers precede, and the right leftovers follow, the output vector (inside $m_{(n)}$). After forming the signed sum that reflects all possible permutations, we go ahead and follow the same procedure for all of the other possible splittings. 

{\bf Example 1.c.} For example, the splitting 
\[ \{ a_1\}\;\;\; \{ a_2\} \;\;\; \{ a_3\} \;\;\; |\;\;\; \{ b_1\}\;\;\;\{ b_2,b_3,b_4\}\;\;\; \emptyset \]
gives rise to the following two terms,
\[ \pm m_{(4)} (a_1,b_1, m_{(1|3)}(a_2|b_2,b_3,b_4),a_3)\]
and
\[  \pm m_{(4)} (b_1,a_1, m_{(1|3)}(a_2|b_2,b_3,b_4),a_3). \]
Note that the inner order of the $a$'s and the $b$'s are preserved. In order to complete the substitution
\[ \{ m_{(4) }\} \{ m_{(1|3)}\}\{ a_1,a_2,a_3 | b_1,b_2,b_3,b_4\} ,\]
the summands coming from each of the remaining five splittings need to be added.

When the type $(n)$ of the leftmost map is replaced by a generic partition $\pi$, things get a bit more complicated. If there is more than one slot, then can we simply substitute the ordered lists constituting $\{ m_{(i_1 | \cdots | i_k)} \} \{ a^{(1)}|\cdots |a^{(k)}\}$ into any set of (consecutive) arguments of $m_{\pi}$ of the same length, then add up? The answer is, not always! We are only allowed to consider those configurations of 
$\{ m_{(i_1 | \cdots | i_k)} \}$ inside $m_{\pi}$ {\it that give rise to the partition} $(n_1|\cdots |n_k)$ of the vectors. If there is {\it no} such configuration, then we have the empty sum, or zero.

{\bf Example 2.} Let us compute a new example,
\[ \{ m_{(2|1)}\} \{ m_{(1|3)}\} \{ a,b|c|d,e,f\} .      \]
Note that the substitution $(\, 2| (1|3)\, )$ gives rise to the partition $(2|1|3)$ of the vectors, whereas $(\, (1|3),1|1\, )$ and $(\, 1,(1|3)|1)$ do not. Then the result is just
\[ \pm m_{(2|1)}(a,b| m_{(1|3)}(c|d,e,f)\, ).\]

The number of terms in a certain composition is easy to compute. Assume that the set of tri-partitions
\[ \{ u_1,i_1,v_1\} ,\dots, \{ u_k,i_k,v_k\}\]
of $n_1,\dots,n_k$ is fixed, where $k\geq 1$ and all numbers are nonnegative integers. First let us find the number of permutations of left (or right) leftovers. 

{\bf Lemma 1} Let $S_1,\dots,S_k$ be ordered sets with $u_1,\dots,u_k$ elements respectively, with $k\geq 1$, and $u_j\geq 0$ for all $1\leq j\leq k$. Then the number of permutations of the elements of $S_1\cup\cdots\cup S_k$ that preserve the inner order of each $S_j$ (called {\it unshuffles}) is given by the product
\[ \left( \begin{array}{c} u_1+\cdots +u_k \\ u_1\end{array}\right)  \left( \begin{array}{c} u_2+\cdots +u_k \\ u_2\end{array}\right)\cdots \left( \begin{array}{c} u_{k-1}+u_k \\ u_{k-1}\end{array}\right) \left( \begin{array}{c} u_k \\ u_k\end{array}\right) .
\]

Once we count all possible values of $u_j$ in the $j$th slot, then the $v_j$'s will fall into place, since the middle value in the tri-partition is $i_j$ and the total is $n_j$. Then 

{\bf Lemma 2} The total number of terms in the substitution
\[ \{ m_{(n)} \} \{ m_{(i_1|\cdots |i_k)}\} \{ a_1^{(1)}\cdots a_{n_1}^{(1)} | \cdots |  a_1^{(k)}\cdots a_{n_k}^{(k)} \} ,\]
where $n=\sum n_j-\sum i_j+1=\sum u_j+\sum v_j+1$, is
\begin{eqnarray} &&\sum_{u_1=0}^{n_1-i_1}\cdots\sum_{u_k=0}^{n_k-i_k}\nonumber \\ &&
\left( \begin{array}{c} u_1+\cdots +u_k \\ u_1\end{array}\right) \left( \begin{array}{c} u_2+\cdots +u_k \\ u_2\end{array}\right)\cdots \left( \begin{array}{c} u_{k-1}+u_k \\ u_{k-1}\end{array}\right) \left( \begin{array}{c} u_k \\ u_k\end{array}\right)\nonumber \\ &&
\cdot \left( \begin{array}{c} v_1+\cdots +v_k \\ v_1\end{array}\right)  \left( \begin{array}{c} v_2+\cdots +v_k \\ v_2\end{array}\right)\cdots \left( \begin{array}{c} v_{k-1}+v_k \\ v_{k-1}\end{array}\right) \left( \begin{array}{c} v_k \\ v_k\end{array}\right) \nonumber \end{eqnarray}
(we are given only the $i_j$'s and $n_j$'s; note that $v_j=n_j-u_j-i_j$ in the formula).
It is not too hard to write a compact formula for the number of terms in case of the most general substitution.

\subsubsection{Symmetries}

Partitions with 0-slots may be used to obtain the symmetries among the generating identities of algebras in the following manner. 

What is a symmetry? In our context, any generating identity that involves terms differing from each other only by permutations of the arguments (vectors) may be called a {\it symmetry}; in our experience, the permutations tend to occur within one map and not a composition of two or more maps (we know one when we see one!). For example, in a commutative associative algebra with bilinear map $m_{(2)}$, the commutativity identity $m_{(2)}(a,b)-m_{(2)}(b,a)=0$ is a symmetry, whereas the associativity identity $m_{(2)}(m_{(2)}(a,b),c)-m_{(2)}(a,m_{(2)}(b,c))=0$ is not, since it  exhibits no change of order in the arguments (also, there are two levels of maps in each term). We claim that this is only a superficial distinction. What if we define, in addition to $m_{(2)}$, a linear map $m_{(1|0)}$ that simply sends a vector to itself and is composed in the funny way indicated? Then
\[ \{ m_{(2)}\} \{ m_{(1|0)}\} \{ a|b\} \]
is evaluated by placing $m_{(1|0)}$ into the first slot of $m_{(2)}$ only (the other possibility gives rise to the partitioning $(2|0)$ of the vectors, instead of $(1|1)$), and {\it not choosing $b$} for substitution into the second slot of $m_{(1|0)}$. We do not have many options as to how to substitute $a$ into the first slot, but in how many ways can we {\it not} substitute $b$ into the second? Let us split the $b$-slot into three sets in all possible manners:
\[ \{ b\}\;\;\;  \emptyset\;\;\; \emptyset,\;\;\;\mbox{or}\;\;\; \emptyset\;\;\;\emptyset\;\;\; \{ b\} .\]
Thus $b$ appears both as a left and a right leftover vector, and the evaluation is given by
\[ m_{(2)}(m_{(1|0)}(a|\emptyset ),b)- m_{(2)}(b,m_{(1|0)}(a|\emptyset ))=m_{(2)}(a,b)-m_{(2)}(b,a)=0.\]
This manifestation of symmetry makes it possible for us to write all-quadratic expressions and use mega identities in a compact way. In particular, (weakly) homotopy Gerstenhaber algebras are distinguished among all homotopy algebras by the inclusion of $m_{(1|0)}$ as the identity map described above. Thus we are able generalize the homotopy associativity of $A_{\infty}$ algebras to homotopy super commutativity and associativity.

\section{Homotopy algebras}

\subsection{Strongly homotopy associative, pre-Lie, and Lie algebras}

\subsubsection{The suspension operator}

Let $m$ be again a formal sum of infinitely many homogeneous maps in the Hochschild space $C^{\bullet}(V,V)$. We want to define another similar map, $\tilde{m}$, by changing each summand by a factor of $\pm 1$. This factor will depend not only on the degree of homogeneity of the summand but also on the arguments of the map.

Following Getzler and Jones \cite{GJ} and mixing in some braces, we describe the construction as follows (also see \cite{Jon}). Let $s$ be the usual {\it suspension operator} (or {\it parity reversion operator})
\[ s:\bigoplus_{j\in \mbox{\bf Z}}V^j\rightarrow \bigoplus_{j\in \mbox{\bf Z}}(sV)^j,\]
where $(sV)^j=V^{j-1}$. Stripped of the grading, $s$ is simply the identity map, hence we have $d(s)=0$ but $|s|=-1$. Similarly,
\[ s^{\otimes k}=\{ s,\dots, s\} \]
sends $V^{k}$ to $(sV)^{k}$, and reduces the super degree by $k$. The tensor product on $sV$ is related to the one on $V$ by the following formula:
\begin{eqnarray} && \{ sa_1,\dots ,sa_k\} \nonumber \\
&=& (-1)^{|s|(|sa_1|+\cdots +|sa_{k-1}|)} \{ i_{(1)},\dots,i_{(1)},s\} \{ sa_1,\dots ,sa_{k-1},a_k\} \nonumber \\
&=& (-1)^{|s|(|sa_1|+\cdots +|sa_{k-1}|)+|s| (|sa_1|+\cdots +|sa_{k-2}|)} \{ i_{(1)},\dots,i_{(1)}, s, s\} \{ sa_1,\dots ,sa_{k-2}, a_{k-1},a_k\} \nonumber \\ && \vdots \nonumber \\
&=& (-1)^{(k-1)|a_1|+(k-2)|a_2|+\cdots +1\cdot |a_{k-1}|+\frac{k(k-1)}{2}} 
\{ s, \dots, s\} \{ a_1,\dots,a_k\} \nonumber \end{eqnarray}

We now define an operator $b_{(k)}$ on $(sV)^{k}$ that makes the diagram
\[ \begin{array}{ccccc} {} & {} & m_{(k)} & {}& {} \\ {} & V^{k} & \rightarrow & V & {}\\ s^{\otimes k} & \downarrow & {} & \downarrow & s \\ {} & (sV)^k & \rightarrow & sV & {}\\
{} & {} & b_{(k)} & {}& {} \end{array} \]
commute. That is, we define $b_{(k)}$ by
\begin{eqnarray} &&  \{ b_{(k)}\}\{ sa_1,\dots,sa_k\} \nonumber \\
&=& \{ s\} \{ m_{(k)}\} \{ s^{-1},\dots,s^{-1}\} \{ sa_1,\dots,sa_k\} \nonumber \\
&=& (-1)^{(k-1)|a_1|+(k-2)|a_2|+\cdots +1\cdot |a_{k-1}|+\frac{k(k-1)}{2}} \{ s\} \{ m_{(k)}\} \{ s^{-1},\dots,s^{-1}\} \{ s, \dots, s\} \{ a_1,\dots,a_k\}\nonumber \\
&=& \{ s\} \left[ (-1)^{(k-1)|a_1|+(k-2)|a_2|+\cdots +1\cdot |a_{k-1}|+\frac{k(k-1)}{2}}
\{m_{(k)}\}  \{ a_1,\dots,a_k\}\right] .\nonumber \end{eqnarray}
We call the expression inside the boxy brackets $\tilde{m_{(k)}}$, and note that we may attribute any properties of $b$, possibly such as $\{ b\}\{ b\} =0$, to the map $\tilde{m}$.

This well-known trick is often employed to make all summands of $m$ odd with respect to the new degree, so that $\tilde{m}$ can be manipulated as one entity. The super degree of $b$, denoted by $||\; ||$ in \cite{Akm2}, is computed from the map's original definition:
\[ ||b_{(k)} ||=|s|+|m_{(k)}|+k|s^{-1}|=|m_{(k)}|+k-1=|m_{(k)}|+d(m_{(k)}).\]
For example, in the definition of $A_{\infty}$ algebras of the next subsection, we have $|m_{(k)}|\equiv k$ (mod 2). Then all $b_{(k)}$, represented by the $\tilde{m}_{(k)}$, become odd operators.

\subsubsection{$A_{\infty}$ algebras}\label{ainf}

We will use the terms ``$A_{\infty}$'' and ``strongly homotopy associative'' interchangeably. Introduced by Stasheff in \cite{Sta1} and \cite{Sta2}, and by now a well-known construct, an $A_{\infty}$ algebra structure on a super graded vector space $V$ is defined by a mega operator $m=m_{(1)}+m_{(2)}+\cdots m_{(n)}+\cdots \in \mbox{\it Hom}(TV,V)$, with $|m_{(k)}|$ having the same parity as $k$, and satisfying 
\[ \{ \tilde{m}\}\{ \tilde{m}\} =0 .\]
This identity is split into sub-identities
\[ \{ \tilde{m}\}\{ \tilde{m}\} \{ a_1,\dots,a_n\} =0\;\;\;\mbox{for all $n\geq 1$,} \]
each of which has finitely many terms, namely,
\[ \sum_{j+k=n+1}\{ \tilde{m_{(j)}}\}\{ \tilde{m_{(k)}}\} \{ a_1,\dots,a_n\} =0 .\]
In \cite{Akm2}, we showed that these are equivalent to
\[ \sum_{j+k=n+1} (-1)^j \{ m_{(j)}\}\{ m_{(k)}\}\{ a_1,\dots,a_n\}  =0.  \]

{\bf Example 1.} \cite{Akm2} It is easily checked using the last equation that given any associative algebra $A$ with bilinear product $m$, the maps $m_{(n)}:A^{\otimes n}\rightarrow A$, defined by zero if $n$ is odd and the unambiguous $(n-1)$-iterated product if $n$ is even, constitute a nontrivial $A_{\infty}$ structure on $A$.

Stasheff's statement \cite{Sta2} (see also \cite{GJ}, \cite{Jon}) that an $A_{\infty}$ structure on $V$ corresponds to a square zero coderivation (of degree $-1$) on $T(sV)$ now follows in one easy step: on both sides, the statement amounts to $\{ \tilde{m}\}\{ \tilde{m}\} =0$ in our notation.

\subsubsection{Pre-$L_{\infty}$ and $L_{\infty}$ algebras}\label{linf}

Strongly homotopy Lie algebras were introduced by Schlessinger and Stasheff in \cite{SS}. Given an $A_{\infty}$ algebra, one may consider the effects of anti-symmetrizing each map $m_{(n)}$. In \cite{Akm2}, we showed that the new structure $\ell =\ell_{(1)}+\ell_{(2)}+\cdots +\ell_{(n)} +\cdots$ satisfies the usual $L_{\infty}$ identities (\cite{LS},\cite{LM}) because
\[ \{ \tilde{m}\}\{ \tilde{m}\} \{ a_1\} \cdots\{ a_n\} =0\]
holds for all $n$. Of course, this particular $L_{\infty}$ construction had been in literature for some time. But then, the new formulation of the $L_{\infty}$ identities led to a simple definition of a pre-$L_{\infty}$ algebra: although a map $m$ need not satisfy the identities
\[ \{ \tilde{m}\}\{ \tilde{m}\} \{ a_1,\dots,a_n\} =0\;\;\;\mbox{for all $n$,} \]
it is enough to find the {\it right pre-$L_{\infty}$} condition
\[ \{ \tilde{m}\}\{ \tilde{m}\} \{ a_1,\dots,\{ a_{n-1}\}\{a_n\}\, \} =0\;\;\;\mbox{for all $n$} \]
in order to conclude that the anti-symmetrization is $L_{\infty}$ (the additional braces around $a_{n-1}$ and $a_n$ indicate, by definition of multibraces, that the signed sum over the permutations of these two elements is taken). For $m=m_{(2)}$ and $n=3$ this is indeed the usual definition of a right pre-Lie algebra. Similarly, one may look for the {\it left pre-$L_{\infty}$} identities
\[ \{ \tilde{m}\}\{ \tilde{m}\} \{\, \{ a_1\}\{ a_2\},\dots,a_n\} =0\;\;\;\mbox{for all $n$} ,\]
or any identity of this type which falls short of $A_{\infty}$ but leads to the $L_{\infty}$ condition above (after being summed over permutations). 

{\bf Example 1.} The anti-symmetrization of the example in Subsection~\ref{ainf} will lead to an $L_{\infty}$ algebra structure on any associative algebra.

{\bf Example 2.} \cite{KVZ, Vor, Akm3} The maps $m_{(1)}$, $m_{(1|1)}$, $m_{(1|1|1)}$,... of the homotopy G-algebra structure described in Part~4 form a pre-$L_{\infty}$ algebra.

{\bf Example 3.} \cite{Akm1, Kra} The higher order differential operators $\Phi_B^1$, $\Phi_B^2$, $\Phi_B^3$,... introduced by the author were shown to constitute an $L_{\infty}$ algebra by Kravchenko in case of a super commutative base algebra. Here $B$ is an odd, linear, square-zero map that is not necessarily a second-order differential operator.

We note that an $L_{\infty}$ algebra structure on a vector space $V$ corresponds to a square-zero codifferential on the exterior coalgebra $\Lambda (sV)$ (see e.g.~\cite{LM}).

\subsection{Generic (weakly) homotopy algebras} \label{homalg}

Among the many circulating definitions of a ``weakly homotopy algebra'' structure on the underlying vector space $V$, we will adopt one that is given by the almost obvious mega identity 
\[ \{ \tilde{m}\}\{ \tilde{m},\tilde{m},\dots \} =0.\]
The mega map $m$ is best thought of as a formal sum of partitioned multilinear maps $m_{\pi}$, one for each partition~$\pi$, some of which may be identically zero. Then the ``algebra of partitions'', introduced in \cite{Akm3} and denoted by ${\cal P}$, determines the sub-identities. For every $\pi$, we find all expressions in ${\cal P}$ of the form $\pi'*[\pi_1,\dots,\pi_k]$, necessarily finitely many, that exhibit $\pi$ as a summand in the final product (partitions have coefficient one or zero in all products of basis elements by definition). These  products correspond to the substitutions of partitioned maps in Subsections~\ref{parmaps} and \ref{parbraces}. Then we sum over the compositions modeled after these products, and set the result equal to zero (some degree restrictions apply). For example, the partition $(1|1)$ occurs as a term only in the partition products $(1)*(1|1)$, $(1|1)*(1)$, $(2)*(1|0)$, and $(2)*(0|1)$, but nowhere else. Therefore, the sub-identity corresponding to the partition $(1|1)$ is written as
\[ \left[  \{ m_{(1)}\} \{ m_{(1|1)}\} \pm \{ m_{(1|1)}\}\{ m_{(1)}\} \pm 
\{ m_{(2)}\} \{ m_{(1|0)}\} \pm \{ m_{(2)}\} \{ m_{(0|1)}\} \right] \{ a|b\} =0.\]

The mega map and the mega identity above were proposed by the author \cite{Akm3} as the ``strongly homotopy Gerstenhaber algebra'' or ``$G_{\infty}$ algebra'' structure on a super graded vector space, following the terminology of the articles \cite{GV}, \cite{KVZ}, and \cite{Vor}. However, we now propose that this identity be known as the generic ``homotopy algebra'' structure, whereas the term ``(weakly) homotopy $G$ algebra'' be reserved for 
homotopy algebras with $m_{(1|0)}$ as the identity map that chooses one element from one slot of vectors and nothing from the very next (there may be several ways of choosin' nuttin') and $m_{(0|1)}$ as the zero map. In particular, the anti-symmetrization of $m_{(1|1)}$ becomes the G-bracket and the bilinear product $m_{(2)}$ becomes the super commutative associative map of the $G$-algebra. The main reason is the simplicity and all-encompassing quality of the more general identity: almost all structures studied under the homotopy umbrella are defined by identities comprising sums of compositions of multilinear maps (there are NO new identities!). Although we are mostly concerned with compositions of one mega map with itself, the necessity of composing such a map with several copies of itself has naturally arisen (again see \cite{GV}, \cite{KVZ}, and \cite{Vor}).

Due to its simplicity, the identity easily lends itself to modifications and generalizations. Substructures are easy to identify (examples will follow), symmetries can be arranged at will, and the partition algebra as well as the depth of the compositions (currently two) may be replaced. 

A note on the grading: we define \cite{Akm3}
\[ ||b_{\pi}||=d(\pi )+\bar{d}(\pi )+|m_{\pi}|+1 ,\]
at least up to congruence modulo two, where $\bar{d}(\pi )$ denotes the number of slots in the partition minus two. In other words, we add the super degree of $|m_{\pi}|$ to the number of arguments and slots of $\pi$, and adopt the resulting parity. The $\bar{d}$-degree will resurface in the study of topological vertex operator algebras. 

There is a generalization of the isomorphism of Subsection~\ref{isom} and the accompanying structures.  The tensor bialgebra $TV$ is replaced by the {\it partitioned tensor bialgebra} $T_{\mbox{\it par}}V=\bigoplus_{\pi}T_{\mbox{\it par}}^{\pi}V$, where every tensor product $a_1\otimes \cdots \otimes a_n$ is now grouped (barred) according to the partition $\pi$. For example, $T_{\mbox{\it par}}^{(1|1)}V$ and $T_{\mbox{\it par}}^{(2)}V$ are different summands of this bialgebra. The bilinear product on $T_{\mbox{\it par}}V$ is defined by the merging of the partitions that mark the basis elements as defined in Section~\ref{parmaps}. Namely,
\[ \{ M\}' \{ \, \{ a^{\pi_1}\} \{ a^{\pi_2}\}\, \} =\{\, \{ a^{\pi_1},a^{\pi_2}\} \,\}' .\]
As for the coproduct $\Delta$, the splitting is now over the sum of all ordered pairs of partitions that ``add up to'', or ``merge into'', the partition $\pi$, instead of $(n)$:
\begin{eqnarray} && \{ \Delta\}' \{ \, \{ a^{\pi}\} \,\}' \nonumber \\
&=& \{ \, \{ i,  i\}\,  \}' \{ \, \{ a^{\pi}\} \,\}' \nonumber \\
 &=&  \sum_{(\pi_1,\pi_2)=\pi} \{ \, \{ i_{\pi_1},i_{\pi_2}   \}\, \}' \{ \, \{ a^{\pi}\} \,\}' \nonumber \\
&=&  \sum_{\{ a^{\pi_1},a^{\pi_2}\} =\{ a^{\pi}\} } \{ \, \{ a^{\pi_1}\} , \{ a^{\pi_2}\}\,\}' .\nonumber \end{eqnarray}
We will not go through the calculations, which remain essentially the same, thanks to the braces notation.

\subsection{Classical and (weakly) homotopy Gerstenhaber algebras}

A {\it (classical) Gerstenhaber algebra} is a super graded vector space $V$ with a super commutative, associative bilinear product (suppressed) and an odd bracket $[\; ,\; ]_G$ satisfying

\noindent (i) $[a,b]_G=-(-1)^{(|a|-1)(|b|-1)}[b,a]_G$ (antisymmetry with suspended grading),

\noindent (ii) $[a,[b,c]_G]_G=[[a,b]_G,c]_G+(-1)^{(|a|-1)(|b|-1)}[b,[a,c]_G]_G$ (Leibniz property with suspended grading), and

\noindent (iii) $[a,bc]_G=[a,b]_Gc+(-1)^{(|a|-1)|b|}b[a,c]_G$ (Poisson rule with mixed grading).

The name comes from the original construction of Gerstenhaber \cite{Ger} on the cohomology of the Hochschild algebra of a super commutative associative algebra. Namely, the super anti-symmetrization of the braces $\{ \, \}\{ \, \}$ mentioned earlier, and the ``dot product'' of multilinear maps, descend to a G-bracket and a commutative associative product with the above properties respectively.

We defined a {\it homotopy Gerstenhaber algebra} in Section~\ref{homalg} as a homotopy algebra with a specific operator $m_{(1|0)}$ that makes the cohomology into a commutative associative algebra and therefore can be thought of as an $A_{\infty}$ algebra with additional structure. There are already such structures in the litearture. For example, in \cite{Mar2}, Markl mentions ``strongly homotopy associative commutative algebras'', a.k.a.~$C_{\infty}$ or balanced $A_{\infty}$ algebras, citing \cite{Mar1} and \cite{Kad}. 

Besides the homotopy G-algebra structure described in Section~\ref{homalg}, another natural way to define a homotopy Gerstenhaber algebra would be via higher brackets $g=\sum g_{\pi}$ (plus some homotopy G-algebra map $m$ in the first sense, to take the place of the super commutative and associative product), one of which directly descends to the G-bracket on the $m_{(1)}$-cohomology. This not only resonates with the way $L_{\infty}$ algebras are defined, but also happens to be the structure on a topological vertex operator algebra (as a homotopy BV-algebra) as we will see. In the second definition, we may want to use the term {\it homotopy Gerstenhaber bracket algebra}, or homotopy GB-algebra, to underline the difference. Depending on example, some general identities resembling those of the classical G-algebras may be imposed on $g$ and $m$. We do not have this problem in the definition of $A_{\infty}$ algebras, for example, because these structures both ``look like'' and ``descend to'' associative algebras.

Two substructures of homotopy G-algebras $(V,m)$ were mentioned in \cite{Akm3}, namely the $A_{\infty}$ subalgebra
\[ m_{(1)}+m_{(2)}+m_{(3)}+\cdots\]
and the pre-$L_{\infty}$ algebra
\[ m_{(1)}+m_{(1|1)}+m_{(1|1|1)}+\cdots .\] For the first, the mega identity is reduced to the usual 
\[ \{ \tilde{m}\}\{ \tilde{m}\} =0 ,\]
and for the second, it is replaced by the right pre-$L_{\infty}$ condition
\[ \{ \tilde{m}\}\{ \tilde{m}\} \{ a_1,\dots,\{ a_{n-1}\}\{a_n\}\, \} =0\;\;\;\mbox{for all $n$}. \]

\subsection{Higher order differential operators} \label{hodo}

The notion of a {\it higher order differential operator} $B\in \mbox{\it Hom}(V,V)$ for a noncommutative, nonassociative algebra $(V,m_{(2)})$ was introduced by the author in \cite{Akm1}. The definition was modeled after Koszul's \cite{Kos} for commutative, associative algebras. We include the definition for reference ($m_{(2)}$ suppressed). Let
\begin{eqnarray} \Phi_B^1(a)&=& B(a)\nonumber \\
\Phi_B^2(a,b)&=&\Phi_B^1(ab)-\Phi_B^1(a)b-(-1)^{|B|\, |a|}aB(b) \nonumber \\
\Phi_B^{r+2}(a_1,\dots,a_r,b,c)&=&\Phi_B^{r+1}(a_1,\dots,a_r,bc)-\Phi_B^{r+1}(a_1,\dots,a_r,b)c
\nonumber \\
&&-(-1)^{|b|(|B|+|a_1|+\cdots +|a_r|)}b\Phi_B^{r+1}(a_1,\dots,a_r,c),\;\;\; r\geq 1. \nonumber \end{eqnarray}
Each bracket measures the previous one's failure to be a derivation of $m_{(2)}$ in the last argument; then we say $B$ is a {\it differential operator of order (at most) $r$} if $\Phi_{B}^{r+1}$ identically vanishes. For a partitioned-map presentation of the formula see \cite{Akm3}. Later, Kravchenko showed in \cite{Kra} that the higher brackets $\Phi_B^r$ form an $L_{\infty}$ algebra structure in case $V$ is super commutative. A generalization of the maps $B$ and $m_{(2)}$ in the $\Phi$-operators to arbitrary multilinear maps is given in~\cite{Akm2}.

\subsection{Classical and weakly homotopy Batalin-Vilkovisky algebras}

The usual definition of a {\it (classical) Batalin-Vilkovisky algebra} (or {\it BV-algebra}) is a Gerstenhaber algebra (with a super commutative associative bilinear product $m_{(2)}$) that is obtained by G-bracketing the bilinear map $m_{(2)}$ with an odd, linear map $\beta$, which is also required to be a second order differential operator with respect to $m_{(2)}$. More precisely, we construct
\[ [a,b]_{\beta}=(-1)^{|a|}[\beta ,m_{(2)}]\{ a,b\} =(-1)^{|a|}\Phi_{\beta}^2(a,b),\]
where the bracket in the middle is the Gerstenhaber bracket of two multilinear maps, and prove that the result is a classical G-algebra (see \cite{Ger, Akm2}).

A {\it weakly homotopy Batalin-Vilkovisky algebra} (or {\it differential BV-algebra}, as it has been called in literature) should be a super graded vector space with a number of higher multilinear maps, in particular an odd, square-zero operator $m_{(1)}$ that is a derivation of some even, bilinear product $m_{(2)}$, and an odd, linear, square-zero operator $B$ that is a second-order differential operator with respect to $m_{(2)}$. We require $m_{(2)}$ be super commutative and associative up to homotopy (although homotopy associativity is needed only to fill out the definition of the classical algebra in the cohomology), and the commutator $[B,m_{(1)}]$ be a diagonalizable operator~$L$ (possibly zero) that commutes with the two degree operators. Under these conditions, it was shown in \cite{Akm1} that all these maps collapse to a classical BV-algebra structure on the $m_{(1)}$-cohomology. 

We may, and indeed will, define a stronger notion of a homotopy BV-algebra based on the two prominent examples to be discussed at the end of the article. We will require a homotopy G-algebra structure and a map $B$ with all the trimmings described above. (Via their mega identity, homotopy G-algebras generalize super commutative associative algebras.) How many classical G-algebra structures does this set-up generate? Well, there is the G-algebra on the $m_{(1)}$-cohomology, coming from the anti-symmetrization of $m_{(1|1)}$, and also the one induced directly by the bracket $[B,m_{(2)}]$ on the same space. This is exactly the situation in a topological vertex operator algebra. In the meantime, we will be freely using the G-bracket notation with no indices on the Hochschild complex to denote the usual anti-symmetrization of simple composition. Too many G's? We also use the name ``Gerstenhaber'' occasionally to refer to a certain mathematician with vast contributions to the field!

{\bf Definition.} A {\it weakly homotopy Batalin-Vilkovisky algebra} is a graded vector space $V$ together with a homotopy G-map $m=\sum_{\pi}m_{\pi}$ and an odd, linear, square-zero, second-order differential operator $B$ with respect to $m_{(2)}$  such that the even linear operator
\[ L=[B,m_{(1)}]\]
is diagonalizable. As in vertex operator algebra theory, we will call the eigenvalues of $L$ {\bf weights}, and note that they produce a third grading on the Hochschild complex of $V$, by definition commuting with the other two. (In many examples, $L$ will be zero.)

{\bf Example 1.} A classical BV-algebra structure on a graded vector space $V$ endows $C^{\bullet}(V,V)$ with a homotopy BV-algebra structure: Subsection~\ref{bvonhoch}.

{\bf Example 2.} A topological vertex operator algebra has a homotopy BV-algebra structure: Subsection~\ref{bvontvoa}.

This definition will give us a {\it higher BV bracket}
\[ [\; ,\; ]_{B}^{\pi}=[B,m_{\pi}]\]
for every partition $\pi$, collectively satisfying identities obtained by G-bracketing those of the homotopy G-algebra structure with $B$. We defined a homotopy GB-algebra previously as a formal sum of brackets satisfying these identities, just as an $L_{\infty}$ algebra structure is a sum of brackets.

{\bf Theorem 1.} A weakly homotopy Batalin-Vilkovisky algebra descends to a classical BV-algebra on $Ker(m_{(1)})/Im(m_{(1)})$. 
 
{\it Proof.} First note that $L$ super commutes with $B$ and $m_{(1)}$. The operator $L$ preserves the kernel of $m_{(1)}$ as a result, and any homogeneous vector $a\in Ker(m_{(1)})$ satisfies
\[ Bm_{(1)}(a)+m_{(1)}Ba=La=\lambda a\]
for some weight $\lambda$. Thus when $\lambda\neq 0$ the vector $a$ falls in $Im(m_{(1)})$ and its image in the cohomology is zero. In other words, $L$ acts with weight zero on the cohomology. Then $B$ and its commutator with $m_{(2)}$ act on the cohomology as well. For the verification of the actual BV identities, see \cite{Akm1}. 

\section{The homotopy algebras of the Hochschild space and topological vertex operator algebra}

\subsection{The homotopy G-algebra structures}

\subsubsection{On the Hochschild complex} \label{homonhoch}

We will review the information in \cite{GV} and \cite{Akm3} about the homotopy G-algebra structure on $C^{\bullet}(V,V)$ for an associative algebra $(V,m_{(2)})$. The nonzero operators are as follows:

\noindent $\{ M_{(1)}\}\{ x\} =[m_{(2)},x]$

\noindent $\{ M_{(2)}\} \{ x,y\} =\pm\{ m_{(2)}\}\{ x,y\}$

\noindent $\{ M_{(1|0)}\} =$ identity map on the left slot (for symmetry)

\noindent $\{ M_{(1|1)}\}\{ x|y\} =\{ x\}\{ y\}$ (shown by Gerstenhaber in \cite{Ger} to be right pre-Lie; see\cite{Akm2})

\noindent $\{ M_{(1|n)}\}\{ x|y_1,\dots,y_n\} =\pm\{ x\}\{ y_1,\dots,y_n\}$ for $n\geq 1$

The following identities were observed by Gerstenhaber and Voronov  in \cite{GV} and rewritten in this form in \cite{Akm3}:

\noindent (i) $M_{(1)}^2=0$

\noindent (ii) $(\{ M_{(1)}\}\{ M_{(2)}\}\pm \{ M_{(2)}\} \{ M_{(1)}\} )(2)=0$

\noindent (iii) $(\{ M_{(2)}\} \{ M_{(2)}\} )(3)=0$

\noindent (iv) $(\{ M_{(1)}\}\{ M_{(1|n+1)}\} \pm \{ M_{(1|n+1)}\}\{ M_{(1)}\}$ 

$\pm \{ M_{(2)}\} \{ M_{(1|n)}\} \pm \{ M_{(1|n)}\} \{ M_{(2)}\} )(1|n+1)=0$ ($n\geq 1$)

\noindent (v) $(\{ M_{(2)}\} \{ M_{(1|n)}\} \pm \{ M_{(1|n)}\} \{ M_{(2)}\}$

$\displaystyle{+\sum_{1\leq i\leq n-1}\pm \{ M_{(2)}\} \{ M_{(1|i)},M_{(1|n-i)}\} )(2|n)=0         }$

\noindent (vi) $(\{ M_{(1)}\}\{ M_{(1|1)}\}\pm\{ M_{(1|1)}\}\{ M_{(1)}\}\pm 
\{ M_{(2)}\}\{ M_{(1|0)}\} )(1|1)=0$

One powerful example of the efficiency of the partitioned-braces language is the newly acquired meaning of the formerly stray identity~(v). This represents the sub-identity corresponding to the partition $(2|n)$, and has terms arising from the compositions of all partitioned maps whose composites contain multilinear maps of type $(2|n)$. Note the substitution of two $M$'s into one as a change from the previous examples. This is in fact why the author decided to try a different mega identity from that of $A_{\infty}$ in \cite{Akm3}. If substitutions of more complicated types arise in examples in the future, we will easily amend the mega-identity, now that we know more of the mechanism behind it.

\subsubsection{On a TVOA}

The homotopy G-algebra structure on a TVOA was revealed by Kimura, Voronov, and Zuckerman in \cite{KVZ, Vor}, albeit without any explicitly constructed maps. (For a complete set of defining properties of a TVOA and some that follow from the definition, see Appendix.) The pictures and the lower identities in \cite{KVZ} were translated into the mega identity $\{ \tilde{m}\}\{ \tilde{m},\tilde{m},\dots\}  =0$ in \cite{Akm3}. The similarities between Lian and Zuckerman's earlier explicit maps and relations \cite{LZ} led the author to seek the connections between the two structures and wonder which, if any, was the fundamental one. 

The operators constructed in this fashion were of all partition types with no zeros and $m_{(1|0)}$ as the sole exception. After Voronov's corrections in \cite{Vor}, it turned out that our generic homotopy algebra mega-identity was an inadequate model for the particular weakly homotopy G-algebra structure on TVOA's. Namely, only the infinitely many identities that correspond to partitions $\pi$ where there is no zero slot, and at most two slots are greater than or equal to two, survived. Unfortunately, such maps plus $m_{(1|0)}$ do not form a G-substructure, solely because of $m_{(1|0)}$ (maps of forbidden partition types appear in valid identities). A compact algebraic description of the whole weakly homotopy G-algebra structure on a TVOA is yet nonexistent.  

In \cite{KVZ}, multibraces denote the equivalent of our partitioned maps. In order to preserve the useful braces notation from earlier work, and to emphasize the importance of the partitions, we have changed the notation of Kimura et al.

\subsection{The homotopy BV-algebra structures}

\subsubsection{On the Hochschild complex} \label{bvonhoch}

We will prove in general that any weakly homotopy algebra structure on $V$, including that of a homotopy BV-algebra, can be lifted naturally to the Hochschild complex, by simplifying Getzler's construction in \cite{Get} for $A_{\infty}$ algebras (see also the interpretation in \cite{Akm2}, p.157).

{\bf Theorem 1.} Let $(V,m)$ be a weakly homotopy algebra (possibly with an additional BV operator $\beta$). Then its structure can be lifted to the Hochschild complex $C^{\bullet}(V,V)$ by defining
\[ \{ M_{\pi}\}'\{ \, \{ x^{\pi}\}\, \}' = \{ \, \{ m_{\pi}\}\{ x^{\pi}\}\, \}' .\]
In case of a BV operator $\beta$ on $V$, we define a BV operator $B$ on $C^{\bullet}(V,V)$ by
\[ \{ B \}'\{ \, \{ x\}\,\}' =\{\,\{ \beta\}\{ x\}\, \}'.\]

{\bf Proof.} With the braces notation, there really isn't much to prove. 
\begin{eqnarray} &&\{ M\}'\{ M,M,M,\dots\}' \{  \, \{ x^{\pi}\}\, \}' \nonumber \\
&=&\sum_{\mbox{product $=\pi$}} \{ M_{\pi'}\}'\{ M_{\pi_1'},M_{\pi_2'},M_{\pi_3'},\dots\}' \{  \, \{ x^{\pi}\}\, \}' \nonumber \\
&=&\sum_{\mbox{product $=\pi$}}   \{ \,\{ m_{\pi'}\}\{ m_{\pi_1'},m_{\pi_2'},m_{\pi_3'},\dots\} \{   x^{\pi}\}\, \}' \nonumber \\
&=& \{\,\{ m\}\{ m,m,m,\dots \}\{   x^{\pi}\}\, \}'\nonumber \\
&=& 0.\nonumber \end{eqnarray}
Clearly, $B$ is odd and square-zero if $\beta$ is. Finally, we can easily show that
\[ \{ \Phi_{B}^3\}'\{\,\{ x,y,z\}\,\}' =\{\,\{ \Phi_{\beta}^3\}\{ x,y,z\}\,\}' =0,\]
that is, $B$ is also a second order differential operator (the first $\phi$ operator is constructed with respect to $M_{(2)}$, and the second one with respect to $m_{(2)}$).

Our construction differs from that of Getzler in that we did not define $M_{(1)}$ by
\[ \{ M_{(1)}\}'\{ \, \{ x\}\, \}'=\{ [m,x]\}', \] 
because it would not have worked: the identity $\{ m\}\{ m\} =0$ (or the $\tilde{m}$ version) is not valid in general, and we would have needed it to prove that $M_{(1)}$ is square-zero. 

\subsubsection{On a TVOA} \label{bvontvoa}

The observation of Lian and Zuckerman \cite{LZ} (also see \cite{PS}) that a homotopy Batalin-Vilkovisky structure exists in the cohomology of a topological vertex operator algebra (see Appendix for full description) predates the discovery of a certain homotopy G-structure on the TVOA itself (\cite{KVZ, Vor, Akm3}). We will put forth strong evidence, with the help of the unifying braces notation, that the former is in fact a consequence of the latter in spite of the chronology. In particular, the homotopy BV-algebra structure on a TVOA, maps, relations, and all, is obtained from the homotopy G-algebra structure by G-bracketing all maps and relations with an odd, square-zero, second-order differential operator $B$.

In \cite{Akm1}, we were able to define higher order differential operators in the most general fashion, and show that ``half'' of the modes $u_n$ for each vertex operator $u(z)$ constituted a sequence of differential operators on $V$ of orders (at most) 1,2,3,..., while the rest satisfied the properties of multiplication operators. Under this novel definition, the mode $b_0$ utilized in \cite{LZ} as the BV operator turned out to be a second-order differential operator with respect to the Wick product, resulting in the vanishing of the higher bracket $\Phi_{b_0}^3(u,v,t)$. We then defined ``generalized '' and ``differential'' BV algebras, on the strength of the existence of second-order differential operators, and gave several examples.

\noindent {\bf Dictionary of maps.} The following table exhibits the multilinear operators listed in  the Lian-Zuckerman article \cite{LZ} versus our interpretation of them with 100\% hindsight.

\[ \begin{array}{|l|l|} \hline \mbox{Lian and Zuckerman call it...} & \mbox{In our language, up to $\pm$...} \\
\hline\hline b_0 & B \\ \hline 
\mbox{$Q$ (BRST differential)} & m_{(1)} \\
\mbox{$\cdot$ (dot product, Wick product)} & m_{(2)} \\
n=n(u,v,t) & m_{(3)} \\ 
m=m(u,v) & m_{(1|1)} \\ \hline
[b_0,Q]=L_0 & [B,m_{(1)}]=L \\
\{ \; ,\; \} \;\;\;\mbox{(BV bracket)} & [B,m_{(2)}]=[\; ,\; ]_{B} \\
\mbox{?} & [B,m_{(3)}] \\ 
m'=m'(u,v) & [B,m_{(1|1)}] \\
n''=n''(u,v,t)\;\;\; \mbox{(appears in Poisson-type identity)} & [m_{(2)},[B,m_{(1|1)}]\, ] \\ \hline \end{array} \]

\noindent For the rest of this section we will mostly write maps and identities in our own symbols.

\noindent {\bf Identities or definitions that were known before Lian-Zuckerman.} 

\noindent {\bf 1.} $[m_{(1)},m_{(2)}]=0$ ($Q$ is a derivation of the Wick, or normal-ordered, product)

\noindent {\bf 2.} $\{ m_{(1)}\}^2 =0$ ($Q^2=0$)

\noindent {\bf 3.} $B^2=0$ ($b_0^2=0$)

\noindent {\bf 4.} $L=[B,m_{(1)}]$ (the commutator of $b_0$ with $Q$ is the (Virasoro) degree operator $L_0$, commuting with all homogeneous maps)

\noindent {\bf Explicit maps/identities found by Lian-Zuckerman.} For simplicity we will not always specify the sign changes due to crossing of symbols. 

\noindent {\bf 5.} $\{ m_{(2)}\}^{sym}=\{ m_{(2)}\}\{ m_{(1|0)}\} =[m_{(1)},m_{(1|1)}]$

\noindent Fundamental homotopy G-identity: Wick is super commutative up to homotopy. First explicit appearance in literature was in \cite{LZ}.

\noindent {\bf 6.} $\{ m_{(2)}\} \{ m_{(2)}\}\pm [m_{(1)},m_{(3)}]=0$

\noindent Fundamental homotopy G-identity: Wick is associative up to homotopy. First explicit appearance in literature again in \cite{LZ}.

\noindent {\bf 7.} $[m_{(1)},[B,m_{(2)}]\, ]=0$ 

\noindent $Q$ is a derivation of the BV-bracket. Why?
\begin{eqnarray} && [m_{(1)},[B,m_{(2)}]\, ]\nonumber \\
&=& [\, [m_{(1)},B],m_{(2)}]\pm [B,[m_{(1)},m_{(2)}]\, ]\nonumber \\
&=& [L,m_{(2)}]\pm [B,0]=0.\nonumber \end{eqnarray}

\noindent {\bf 8.} $[\; ,\; ]_{B}^{sym}=[B,m_{(2)}]^{sym}=[m_{(1)},[B,m_{(1|1)}]\, ]$

\noindent The BV-bracket is symmetric up to homotopy. Can be obtained as follows:
\begin{eqnarray} && [m_{(1)},[B,m_{(1|1)}]\, ] \nonumber \\
&=& [\, [m_{(1)},B],m_{(1|1)}]\pm [B,[m_{(1)},m_{(1|1)}]\, ]\nonumber \\
&=& [L,m_{(1|1)}]\pm [B,m_{(2)}^{sym}]\nonumber \\
&=& \pm [B,m_{(2)}]^{sym}.\nonumber \end{eqnarray}

\noindent {\bf 9.} $[u,[v,t]_{B}]_{B}=[\, [u,v]_{B},t]_{B}\pm [v,[u,t]_{B}\,]_{B}$

\noindent Leibniz property. In \cite{Akm2} (pp. 154-155), we showed that this is equivalent to the vanishing of the operator 
\[ \Phi_{B^2}^{3}-[B,\Phi_B^3].\]
Since $B$ is a square-zero, second-order differential operator, we are done.

\noindent {\bf 10.} $[u,\{ m_{(2)}\}\{ v,t\} ]_{B}=\{ m_{(2)}\} \{ [u,v]_{B},t\} \pm \{ m_{(2)}\} \{ v,[u,t]_{B} \}$

\noindent Poisson property for left BV-bracketing with respect to the Wick product. Again in \cite{Akm2} (p.155) we showed that this identity is exactly the condition for the operator $\Phi_B^3$ to be identically zero.

\noindent {\bf 11.} $[\{ m_{(2)}\}\{ u,v\} ,t]_{B}=\{ m_{(2)}\} \{ u,[v,t]_{B}   \} \pm
\{ m_{(2)}\} \{ [u,t]_{B}, v\}$ 

\noindent $\pm [m_{(1)},[m_{(2)},[B,m_{(1|1)}]\, ]\, ]\{ u,v,t\}$

\noindent Poisson property for right BV-bracketing up to homotopy. We do not have a pre-cooked identity that corresponds to this one (due to the lopsided way we have defined our $\Phi$-operators).

\noindent {\bf Two more related identities.}

\noindent {\bf 12.} $\{ m_{(2)}\}\{ \{ u\}\{ v\} ,t\} =0$

\noindent Found in \cite{Akm3}: Wick is left pre-Lie, hence its G-bracket is Lie. The Lie property of the G-bracket was also independently observed by Dong, Li, and Mason in \cite{DLM}.

\noindent {\bf 13.} $[B,\{ m_{(2)}\}\{ m_{(2)}\} ]=[B,[m_{(1)},m_{(3)}]]=[m_{(1)},[B,m_{(3)}]]$

\noindent This identity is obtained by B-bracketing the homotopy associativity statement, and could easily have been discovered by Lian and Zuckerman. The operator $[B,m_{(3)}]$ was similarly within their scope.

\noindent {\bf One more homotopy BV structure.} By the Theorem in Subsection~\ref{bvonhoch}, there is also a homotopy BV-algebra structure with the Hochschild complex of any TVOA as the base space. 

\subsubsection{Explicit construction?}

How do we construct an {\it explicit} homotopy G-map $m$ on a generic TVOA, now that the mystery behind its homotopy BV structure is explained? The key may be found in Lian and Zuckerman's construction of the ``$BV_{\infty}$'' maps in \cite{LZ}, using the relation 
\[ Qb_{-1}+b_{-1}Q=L_{-1}\]
several times. We will show alternative brute-force constructions of the same maps based on the VOA properties displayed in the Appendix and especially utilizing the important identity above. This work is yet unfinished, but completion may become technically possible, due to the existence of an induction argument. The induction is on the $d+\bar{d}$-degree of the multilinear maps we are after (that is, the sum of the numbers of arguments and slots minus three). This particular grading is not preserved under composition; in fact, $d$ is preserved but $\bar{d}$ is decreased by one. With the exception of $(1)$ every partition has a nonnegative $d+\bar{d}$-degree, and the lower maps are already known (see the nontrivial ones below). Now we always have the terms 
\[ \{ m_{(1)}\}\{ m_{\pi'}\}\pm\{ m_{\pi'}\}\{ m_{(1)}\}\]
in the sub-identity corresponding to type $\pi'$. Since $(d+\bar{d})(m_{(1)})=-1$, we observe that the $d+\bar{d}$-degree of $ m_{\pi'}$ (say $n$) is one more than the sum of the two $d+\bar{d}$-degrees in the remaining terms (total $n-1$) of the sub-identity, or at least one more than the $d+\bar{d}$-degree of every remaining map (at most $n-1$). That is, all the other partitioned maps that appear will have been constructed prior to our $m_{\pi'}$! Then we start with the terms in the rest of the identity and try to put their sum into the form $\{ m_{(1)}\}\{ m_{\pi'}\}\pm\{ m_{\pi'}\}\{ m_{(1)}\}$.

First, let us investigate the ``G-bracket'' map $m_{(1|1)}$. It appears in the sub-identity
\[ \{\, \{ m_{(1)}\}\{ m_{(1|1)}\}\pm \{ m_{(1|1)}\}\{ m_{(1)}\}\pm\{ m_{(2)}\}\{ m_{(1|0)}\}\,\} \{ u|v\} =0.\]
We then start with the last term and note that, after cancellations, the expression is in the image of the ``differentiation'' operator $L_{-1}=Qb_{-1}+b_{-1}Q$:
\begin{eqnarray} &&\{ m_{(2)}\}\{ m_{(1|0)}\}\{ u|v\} \nonumber \\
&=& u_{-1}v-v_{-1}u \nonumber \\
&=& [u_{-1},v_{-1}]\mbox{\bf 1} \nonumber \\
&=& \sum_{i\geq 0}(-1)^i(u_iv)_{-2-i}\mbox{\bf 1} \nonumber \\
&=& \sum_{i\geq 0}(-1)^i\frac{1}{-1-i}[L_{-1},(u_iv)_{-1-i}]\mbox{\bf 1} \nonumber \\
&=& \sum_{i\geq 0}\frac{(-1)^{i+1}}{i+1}L_{-1}(u_iv)_{-1-i}\mbox{\bf 1} \nonumber \\
&=& \sum_{i\geq 0}\frac{(-1)^{i+1}}{i+1}(Qb_{-1}+b_{-1}Q)(u_iv)_{-1-i}\mbox{\bf 1} \nonumber \\
&=& Q\left( \sum_{i\geq 0}\frac{(-1)^{i+1}}{i+1}b_{-1}(u_iv)_{-1-i}\mbox{\bf 1}\right) \nonumber \\
&& +\sum_{i\geq 0}\frac{(-1)^{i+1}}{i+1}b_{-1}(\, (Qu)_iv)_{-1-i}\mbox{\bf 1} \pm \sum_{i\geq 0}\frac{(-1)^{i+1}}{i+1}b_{-1}(u_i(Qv)\, )_{-1-i}\mbox{\bf 1} \nonumber \\
&=& [Q, m_{(1|1)}]\{ u|v\} ,\nonumber \end{eqnarray}
where
\[ \{ m_{(1|1)}\}\{ u|v\} =\sum_{i\geq 0}\frac{(-1)^{i+1}}{i+1}b_{-1}(u_iv)_{-1-i}\mbox{\bf 1}
=\sum_{i\geq 0}\frac{(-1)^{i+1}}{i+1}(b_{-1}u_iv)_{-1-i}\mbox{\bf 1} .\]

Similarly, the partition $(3)$ gives rise to the identity
\[ \{\, \{ m_{(1)}\}\{ m_{(3)}\} \pm \{ m_{(3)}\} \{ m_{(1)}\} \pm \{ m_{(2)}\}\{ m_{(2)}\}\,\}
\{ u,v,w\} =0.\]
We start again from the term that has already been defined, i.e.
\begin{eqnarray} && \{ m_{(2)}\}\{ m_{(2)}\}\{ u,v,w\} \nonumber \\
&=&(u_{-1}v)_{-1}w-u_{-1}(v_{-1}w) \nonumber \\
&=& \left( \sum_{i\geq 0}u_{-1+i}v_{-1-i}+v_{-2-i}u_i\right) w-u_{-1}v_{-1}w \nonumber \\
&=& \sum_{i\geq 1}u_{-1+i}v_{-1-i}w+\sum_{i\geq 0}v_{-2-i}u_iw \nonumber \\
&=& \sum_{i\geq 0}u_iv_{-2-i}w+\sum_{i\geq 0}v_{-2-i}u_iw \nonumber \\
&=& \sum_{i\geq 0}\frac{1}{i+1}(L_{-1}u)_{i+1}v_{-i-2}w+\sum_{i\geq 0}\frac{1}{-i-1}(L_{-1}v)_{-i-1}u_iw \nonumber \\
&=& \sum_{i\geq 0}\frac{1}{i+1}(\, (Qb_{-1}+b_{-1}Q)u)_{i+1}v_{-i-2}w+\sum_{i\geq 0}\frac{1}{-i-1}(\, (Qb_{-1}+b_{-1}Q)v)_{-i-1}u_iw \nonumber \\
&=& \sum_{i\geq 0}\frac{1}{i+1}[Q,(b_{-1}u)_{i+1}]v_{-i-2}w+\sum_{i\geq 0}\frac{1}{i+1}(b_{-1}Qu)_{i+1}v_{-i-2}w \nonumber \\
&& -\sum_{i\geq 0}\frac{1}{i+1}[Q,(b_{-1}v)_{-i-1}]u_iv-\sum_{i\geq 0}\frac{1}{i+1}(b_{-1}Qv)_{-i-1}u_iw \nonumber \\
&=& \sum_{i\geq 0}\frac{1}{i+1}Q(b_{-1}u)_{i+1}v_{-i-2}w\pm \sum_{i\geq 0}\frac{1}{i+1}(b_{-1}u)_{i+1}Qv_{-i-2}w \nonumber \\
&& \pm \sum_{i\geq 0}\frac{1}{i+1}(b_{-1}Qu)_{i+1}v_{-i-2}w 
-\sum_{i\geq 0}\frac{1}{i+1}Q(b_{-1}v)_{-i-1}u_iw \nonumber \\ 
&& \pm \sum_{i\geq 0}\frac{1}{i+1}(b_{-1}v)_{-i-1}Qu_iw\pm \sum_{i\geq 0}\frac{1}{i+1}(b_{-1}Qv)_{-i-1}u_iw \nonumber \\
&=& \{ Q\} \{ m_{(3)}\} \{ u,v,w\} \pm \{ m_{(3)}\} \{ Q\} \{ u,v,w\} ,\nonumber \end{eqnarray}
with
\[ \{ m_{(3)}\} \{ u,v,w\} =\sum_{i\geq 0}\frac{1}{i+1}(b_{-1}u)_{i+1}v_{-i-2}w-\sum_{i\geq 0}\frac{1}{i+1}(b_{-1}v)_{-i-1}u_iw.\]

In order to extend the construction to $m$, it remains to show that every ``computed expression'' that we start with is in fact in the image of $L_{-1}$. Or, we need to ascertain that every such expression starts with some mode $u_n$ with $n\neq -1$, as we have
\[ u_n=\frac{1}{n+1}[L_{-1},u_{n+1}].\]

\section{Summary of classification}

We have defined the most general (quadratic) {\bf weakly homotopy algebra} on a super graded vector space $V$ by the relation
\[ \{ \tilde{m}\}\{ \tilde{m},\tilde{m},\dots\} =0\]
where 
\[ m=\sum_{\pi}m_{\pi}\in C_{par}^{\bullet}(V,V),\]
and there are no conditions on the partitioned maps $m_{\pi}$ except possibly for the super degrees. In order to obtain a super commutative, associative algebra in the $m_{(1)}$-cohomology, we impose the restriction
\[ m_{(1|0)}=\mbox{identity and $m_{(0|1)}=$zero}\]
on the homotopy algebra and obtain a {\bf weakly homotopy Gerstenhaber algebra}. Most explicit examples discussed in this article fall into this class. If there is a homotopy G-algebra structure on $V$ as well as an odd, linear, square-zero, second-order differential operator $B$ (with respect to the bilinear map $m_{(2)}$), then we call the resulting form a {\bf weakly homotopy Batalin-Vilkovisky algebra}. But since $G$-bracketing $m$ with $B$ yields a sum $g=\sum [\; ,\; ]_B^{\pi}$ of higher BV brackets on $V$, descending to a classical BV algebra bracket on the cohomology, we also introduce the concept of a {\bf weakly homotopy Gerstenhaber bracket algebra} having a homotopy G-map $m$ and a bracket map $g$. The classical G-algebra structure that exists on the cohomology in both cases is induced by the anti-symmetrization of $m_{(1|1)}$ and the product $m_{(2)}$, but in the second case, we also have the G-bracket $g_{(2)}$ or the BV-algebra bracket $[B,m_{(2)}]$ that descends to the cohomology, together with $m_{(2)}$.

A {\bf strongly homotopy associative algebra} ($A_{\infty}$ algebra) is a homotopy algebra with only nonzero maps of partition type $(n)$, satisfying
\[ \{ \tilde{m}\}\{ \tilde{m}\} =0,\]
and its anti-symmetrization leads to an $L_{\infty}$ algebra, or {\bf strongly homotopy Lie algebra}, of brackets. In fact, any algebra $(V,m)$ that produces $L_{\infty}$ under anti-symmetrization may be called {\bf pre-$L_{\infty}$} and may possibly satisfy 
\[ \{ \tilde{m}\}\{ \tilde{m}\}\{ a_1,\dots,a_{n-2},\{ a_{n-1}\}\{ a_n\}\, \} =0\;\;\; \mbox{or} \;\;\; \{ \tilde{m}\}\{ \tilde{m}\}\{ \,\{ a_1\}\{ a_2\},\dots, a_n \} =0 \]
(the right and left pre-$L_{\infty}$ identities respectively). An $L_{\infty}$ algebra obtained by anti-symmetrization of some $m$ satisfies the relation
\[  \{ \tilde{m}\}\{ \tilde{m}\}\{ a_1\}\dots\{ a_n\} =0.\]
We hope to bring about a new wave of discussions of names and definitions, even though our nomenclature may not be accepted as is. 
\vspace{.3in}

\noindent {\bf Acknowledgments.} I am indebted to Martin Markl, Jim Stasheff, and Sasha Voronov for numerous helpful comments and corrections. Any remaining errors are mine, as the saying goes.


\section{Appendix: Properties of VOA's and TVOA's}

\noindent A {\bf vertex operator algebra (VOA)} is a {\bf Z}-weighted (possibly with additional super {\bf Z}-grading, denoted by superscripts; not assumed here) module
\[ V=\bigoplus_nV_n\]
of the Virasoro algebra that satisfies the following properties:

\noindent {\bf P1.} $V_n=0$ for $n<<0$.

\noindent {\bf P2.} There exists a linear map
\[ V\rightarrow End(V)[\, [z,z^{-1}]\, ],\;\;\; v\mapsto v(z)=\sum v_nz^{-n-1},\]
where $z$ is a formal variable, $v$ is called a {\it state}, and $v(z)$ is called a {\it vertex operator} (or a {\it field}). The linear map $v_n$ is called a {\it mode} of the vertex operator $v(z)$. We will sometimes denote the product $v_nw$ by $v\times_nw$.

\noindent {\bf P3.} For all $v,w\in V$, we have $v_nw=0$ for $n>>0$.

\noindent {\bf P4.} There exists an element {\bf 1}$\in V_0$ such that
\[ \mbox{\bf 1}(z)=\mbox{id}\cdot z^0.\]

\noindent {\bf P5.} There exists $\omega\in V_2$ such that 
\[ \omega(z)=\sum L_nz^{-n-2},\]
where $L_n$ denotes both an element of $Vir$ and the corresponding mode in $End(V)$.

\noindent {\bf P6.} For all $v\in V$, we have
\[ v(z)\mbox{\bf 1}\in V[\, [z]\, ]\]
and
\[ \lim_{z\rightarrow 0}v(z)\mbox{\bf 1}=v.\]

\noindent {\bf P7.} The operator $L_0$ is diagonalizable, with
\[ L_0v=(wt\, v)v\;\;\;\mbox{for all $v\in V_n$}\]
(the eigenvalues are called {\it weights}).

\noindent {\bf P8.} We have
\[ \frac{d}{dz}v(z)=(L_{-1}v)(z)\;\;\;\mbox{for all $v\in V$.}\]

\noindent {\bf P9.} Bracketing with $L_0$ results in
\[ [L_0,v(z)]=(L_0v)(z)+z(L_{-1}v)(z)\;\;\;\mbox{for all $v\in V$.}\]

\noindent {\bf P10.} For all $v,w\in V$, there exists an integer $t>>0$ such that
\[ [v(z_1),w(z_2)](z_1-z_2)^t=0.\]

\noindent Some additional properties that follow from the definition are:

\noindent $\star$ $\displaystyle{(u_mv)_n=\sum_{i\geq 0}(-1)^i\left( \begin{array}{c} m \\ i \end{array} \right) [u_{m-i}v_{n+i}-(-1)^mv_{m+n-i}u_i]}$

\noindent $\star$ $\displaystyle{[u_m,v_n]=\sum_{i\geq 0}(-1)^i(u_iv)_{m+n-i}}$

\noindent $\star$ $\displaystyle{\left( \begin{array}{c} 0 \\ i \end{array} \right) =\left\{
\begin{array}{ll} 1 & \mbox{if $i=0$} \\ 0 & \mbox{if $i\geq 1$} \end{array} \right. }$

\noindent $\star$ $\displaystyle{\left( \begin{array}{c} -1 \\ i \end{array} \right) =(-1)^i\;\;\;\mbox{($i\geq 0$)}}$

\noindent $\star$ $\displaystyle{\left( \begin{array}{c} m \\ i \end{array} \right) =\left\{
\begin{array}{l} 1 \;\;\; \mbox{if $i=0$} \\ m(m-1)(m-2)\cdots (m-i+1)/i!\;\;\; \mbox{if $i\geq 1$} \end{array} \right. }$

\noindent for $m\in${\bf Z}

\noindent $\star$ $L_0-(n+1)\cdot\mbox{id}$ is a derivation of $\times_n$ ($L_0$ is a derivation of Wick)

\noindent $\star$ $[L_0,u_n]=(wt \, u-n-1)\cdot u_n$

\noindent $\star$ $\displaystyle{(u_0v)(z)=[u_0,v(z)]=\sum [u_0,v_n]z^{-n-1}}$

\noindent $\star$ $\displaystyle{\frac{d}{dz}\sum u_nz^{-n-1}=(L_{-1}u)(z)=\sum [L_{-1},u_n]z^{-n-1}=-\sum nu_{n-1}z^{-n-1}}$

\noindent ($L_{-1}$ is also a derivation of Wick)

\noindent $\star$ The mode $u_n$ is a differential operator of order $\leq n+1$ with respect to Wick for $n\geq 0$

\noindent $\star$ Wick is left pre-Lie, or $u_{-1}v_{-1}-v_{-1}u_{-1}-(u_{-1}v)_{-1}+(v_{-1}u)_{-1}=0$

\noindent $\star$ $(u_m\mbox{\bf 1})_n$ has at most two terms:

\noindent $\displaystyle{(-1)^{-n-1} \left( \begin{array}{c} m \\ -n-1 \end{array} \right) u_{m+n+1}}$ if $n\leq -1$ and $\displaystyle{(-1)^n\left( \begin{array}{c} m \\ m+n+1 \end{array} \right) u_{m+n+1}}$ if $m+n\geq -1$

\noindent $\star$ For a VOSA with nonnegative weights, $V_0$ is a (super) commutative, associative algebra under Wick

\noindent $\star$ For $wt\, u =m$ and $wt \, v=n$, all expressions $u_tv$ are zero for $t\geq m+n-$minimal weight

\noindent $\star$ A {\it residue (charge)} $u_0$ is a derivation of all binary products $\times_n$:

\noindent $[u_0,\times_n](v,w)=u_0(v_nw)\pm (u_0v)_nw\pm v_n(u_0w)=0$

\noindent $\star$ Again if $u_0$ is a residue, we have 

\noindent $u_0v_n\mbox{\bf 1}=[u_0,v_n]\mbox{\bf 1} =(u_0v)_n\mbox{\bf 1}$ and $u_0v_mw_n\mbox{\bf 1}=(u_0v)_mw_n\mbox{\bf 1}\pm v_m(u_0w)_n\mbox{\bf 1}$

\noindent $\star$ $\displaystyle{u_n=\frac{1}{n+1}[L_{-1},u_{n+1}]}$ for $n\neq -1$.

\noindent A {\bf topological vertex operator algebra (TVOA)} is a super graded VOA with the following additional properties:

\noindent {\bf P11.} There exists an even vertex operator $F(z)\in V_1$ whose residue $F_0$ is the {\it fermion (ghost) number operator}, or the super degree operator

\noindent {\bf P12.} There exists a weight-one primary (Virasoro-singular) vertex operator $J(z)$ with fermion number one and a square-zero residue $Q=J_0$

\noindent {\bf P13.} There exists a weight-two primary operator $b(z)=\sum b_nz^{-n-2}$ with fermion number $-1$, satisfying
\[ [Q,b(z)]=\omega (z)\]


\begin{thebibliography}{99}

\bibitem{Akm1} F.~Akman, On some generalizations of Batalin-Vilkovisky algebras, {\it J.~Pure Appl.~Algebra} {\bf 120} (1997), 105-141. 

\bibitem{Akm2} F.~Akman, Multibraces on the Hochschild space, {\it J.~Pure Appl.~Algebra} {\bf 167} (2002), 129-163. 

\bibitem{Akm3} F.~Akman, A master identity for homotopy Gerstenhaber algebras, {\it Commun.~Math.~Phys.} {\bf 209} (2000), 51-76.

\bibitem{DL} C.~Dong and J.~Lepowsky, {\it Generalized Vertex Algebras and Relative Vertex Operators}, {\it Progress in Mathematics} {\bf 112}, Boston, Birkh\"{a}user, 1993. 

\bibitem{DLM} C.~Dong, H.~Li, and G.~Mason, Vertex Lie algebras, vertex Poisson algebras, and vertex operator algebras, e-print arXiv:math.QA/0102127.

\bibitem{Ger} M.~Gerstenhaber, The cohomology structure of an associative ring, {\it Ann.~Math.}    {\bf 78} (1963), 267-288.

\bibitem{GV}  M.~Gerstenhaber and A.~A.~Voronov, Higher order operations on Hochschild complex, {\it Functional Anal.~Appl.} {\bf 29 (1)} (1995), 1-6.

\bibitem{Get} E.~Getzler, Cartan homotopy formulas and the Gauss-Manin connection in cyclic homology, {\it Israel Mathematical Conference Proceedings} {\bf 7} (1993), 65-78.

\bibitem{GJ} E.~Getzler and J.~D.~S.~Jones, $A_{\infty}$-algebras and the cyclic bar complex, {\it Ill.~Jour.~Math.} {\bf 34} (1989), 256-283.

\bibitem{HZ} Y.-Z.~Huang and W.~Zhao, Semi-infinite forms and topological vertex operator algebras, {\it Comm. Contemp. Math.} {\bf 2} (2000), 191-241.

\bibitem{Jon} J.~D.~S.~Jones, Lectures on operads. In {\it Quantization, Poisson Brackets and Beyond} (T.~Voronov ed.), {\it Contemporary Mathematics} {\bf 315} (2002), 89-130.

\bibitem{Kad} T.~V.~Kadeishvili, O kategorii differentialnych koalgebr i kategorii $A(\infty)$-algebr, {\it Trudy Tbiliss. Mat. Inst. Razmadze Akad. Nauk Gruzin. SSR} {\bf 77} (1985), 50-70, in Russian.

\bibitem{KVZ} T.~Kimura, A.~A.~Voronov, and G.~J.~Zuckerman, Homotopy Gerstenhaber algebras and topological field theory. In: J.-L.~Loday, J.~Stasheff, and A.~A.~Voronov (eds.), {\it Operads: Proceedings of Renaissance Conferences}, {\it Contemp.~Math.} {\bf 202} (1996), 305-333.

\bibitem{Kos} J.-L.~Koszul, Crochet de Schouten-Nijenhuis et cohomologie, {\it Ast\'{e}risque} (1985), 257-271.

\bibitem{Kra} O.~Kravchenko, Deformations of Batalin-Vilkovisky algebras, Banach Center Publications, Institute of Mathematics, Polish Academy of Sciences, Warszawa, 1999 (e-print arXiv:math.QA/9903191).

\bibitem{LM} T.~Lada and M.~Markl, Strongly homotopy Lie algebras, {\it Comm. Algebra} {\bf 23(6)} (1995), 2147-2161.

\bibitem{LS} T.~Lada and J.~Stasheff, Introduction to SH Lie algebras for physicists, {\it Internat.~J.~Theoret. Phys.} {\bf 32(7)} (1993), 1087-1103.

\bibitem{LZ} B.~H.~Lian and G.~J.~Zuckerman, New perspectives on the BRST-algebraic structure of string theory, {\it Commun.~Math.~Phys.} {\bf 154} (1993), 613-646.

\bibitem{Mar1} M.~Markl, A cohomology theory for $A(m)$-algebras and applications, {\it J.~Pure Appl.~Algebra} {\bf 83} (1992), 141-175.

\bibitem{Mar2} M.~Markl, Homotopy algebras are homotopy algebras, e-print arXiv:math.AT/9907138.

\bibitem{PS} M.~Penkava and A.~Schwarz, On some algebraic structures arising in string theory. In: {\it Perspectives in Mathematical Physics, Conf.~Proc.}, {\it Lecture Notes in Math.~Phys.} {\bf III}, Cambridge, MA, International Press, 1994.

\bibitem{SS} M.~Schlessinger and J.~Stasheff, The Lie algebra structure of tangent cohomology and deformation theory, {\it J.~Pure Appl.~Algebra} {\bf 38} (1985), 313-322.

\bibitem{Sta1} J.~D.~Stasheff, Homotopy associativity of H-spaces, I, {\it Trans. Amer. Math. Soc.} {\bf 108} (1963), 275-292.

\bibitem{Sta2} J.~D.~Stasheff, Homotopy associativity of H-spaces, II, {\it Trans. Amer. Math. Soc.} {\bf 108} (1963), 293-312.

\bibitem{TT} D.~Tamarkin and B.~Tsygan, Noncommutative differential calculus, homotopy BV algebras, and formality conjectures, arXiv.math.KT/0002116.

\bibitem{Vor} A.~A.~Voronov, Homotopy Gerstenhaber algebras, {\it Conference Moshe Flato 1999} (G. Dito and D. Sternheimer, eds.), vol. 2, Kluwer Academic Publishers, the Netherlands, 2000, 307-331.


\end{thebibliography}
\end{document}